# CHARACTERISTIC P GALOIS REPRESENTATIONS THAT ARE PRODUCED BY DRINFELD

Nigel Boston and David T. Ose

University of Illinois at Urbana-Champaign and Bucknell University

June 1997

## 0. Introduction

There are well-known methods of producing representations of the absolute Galois group of a number field. These include the use of elliptic curves, modular forms, and most generally étale cohomology groups of varieties [FM]. There are many conjectures as to which Galois representations are produced this way. For instance, Serre's conjecture [S] states that every odd, irreducible representation of the form $\mathrm{Gal}(\overline{\mathbf{Q}}/\mathbf{Q}) \to GL_2(\overline{\mathbf{F}_p})$ should be associated to a modular form of a particular kind. Here *odd* means that complex conjugation maps to a matrix of determinant $-1$.

In this paper, we consider the analogous case of representations of the absolute Galois group of a field of nonzero characteristic. Suppose that $K$ has characteristic $p \neq 0$. We describe a method, due to Drinfeld, of obtaining representations of the form $\mathrm{Gal}(\overline{K}/K) \to GL_r(\overline{\mathbf{F}_p})$ and address the problem of which representations arise this way. We obtain a fairly complete answer in the case $r = 1$. This has applications to classifying cyclic extensions of $K$ of degree $m$, even when the $m$th roots of unity are not all in $K$. The question of what representations of the form $\mathrm{Gal}(\overline{K}/K) \to GL_r(R)$ ($R$ a discrete valuation ring of equal characteristic with finite residue field) are produced by extending the method of Drinfeld, is addressed in the second author's University of Illinois Ph.D. thesis [O].

## 1. Drinfeld Representations

Let $K$ be a field of characteristic $p$. Suppose that $K$ contains $\mathbf{F}_q$. Define the *Ore ring* to be the set of polynomials in $F$ over $K$, $K\{F\} = \{\sum a_i F^i : a_i \in K\}$, with the noncommutative multiplication $Fa = a^q F$. This ring is also known as the ring of $\mathbf{F}_q$-linear polynomials or alternatively $\mathrm{End}_{\mathbf{F}_q}(\mathbf{G}_a/K)$, with $F$ interpreted as

The first author was partially supported by NSF grant DMS 96-22590, the Sloan Foundation, and the Rosenbaum Foundation. The authors thank God for leading them to results. They also thank Bucknell University for supporting their collaboration, the University of Illinois for enabling them to visit England, the Newton Institute for its hospitality, and Hendrik Lenstra for useful comments on this work.

the Frobenius morphism that sends $x$ to $x^q$. For its basic properties, see Chapter 1 of [G]. Let $g(F) \in K\{F\}$ be of degree $r > 0$. Let $\phi \in \mathbf{F}_q[T]$ be irreducible and of degree $d > 0$. We make the assumption that $\phi(b) \neq 0$, where $b$ is the constant term of $g$. The set $V = \{x \in \overline{K} : (\phi(g(F)))x = 0\}$ (where $Fx = x^q, F^2x = x^{q^2}, ...$) is a vector space over $\mathbf{F}_{q^d}$ of dimension $r$, sometimes called the $\phi$-division points, on which $G_K = \mathrm{Gal}(\overline{K}/K)$ acts (the assumption on $\phi(b)$ ensuring that $(\phi(g(F)))x$ is separable so that $V$ has the claimed cardinality). We therefore obtain a representation $\rho : G_K \to GL_r(\mathbf{F}_{q^d})$. The question we wish to address is what representations arise this way. Such representations will be called *Drinfeld* (but note that Drinfeld modules may be more general).

## 2. A Useful Lemma

Let $g(F) = aF + b$ and so $r = 1$. Then $\rho$ maps to $\mathbf{F}_{q^d}^*$, and hence factors through $\mathrm{Gal}(L/K)$, where $L/K$ is a cyclic extension of degree dividing $q^d - 1$ ($L = K(V)$ in the notation of the introduction - we will denote it by $L_{a,b,\phi}$ in later work). Let $\zeta$ be a root of $\phi$, $K' = K(\zeta)$, and $L' = L(\zeta)$. The extension $L'/K'$ is a Kummer extension since $K' = K(\zeta)$ contains $\mathbf{F}_q(\zeta) = \mathbf{F}_{q^d}$. Thus, $L' = K'(v)$ where $v^{q^d-1} \in K'$, say $v^{q^d-1} = c$.

What we need to know is the following. What is $c$ in terms of $a, b, \zeta$?

**Lemma** With the set-up as above, $c = \frac{(\zeta - b)(\zeta - b^q)...(\zeta - b^{q^{d-1}})}{a^{1+q+...+q^{d-1}}}$.

*Proof* Let $\phi(T) = (T - \zeta)\psi(T)$, so $\psi(T)$ is a polynomial over $\mathbf{F}_q(\zeta) = \mathbf{F}_{q^d}$ of degree $d - 1$. Let $x \neq 0$ satisfy $(\phi(aF + b))x = 0$, so that $L = K(x)$ (since $L/K$ is cyclic) and $L' = K'(x)$.

We claim that if $v = (\psi(aF + b))x$, then $L' = K'(v)$, and most importantly $v^{q^d-1} = (\zeta - b)(\zeta - b^q)...(\zeta - b^{q^{d-1}})/a^{1+q+...+q^{d-1}}$.

This follows from the following identity in $K'\{F\}$. (Here $[q]_k$ stands for $(q^k - 1)/(q - 1)$ and $c_i$ stands for $\zeta - b^{q^i}$.)

$$(a^{[q]_d}F^d - c_0c_1c_2...c_{d-1})\psi(aF + b) = h(F)\phi(aF + b),$$

where

$$h(F) = a^{[q]_{d-1}}F^{d-1} + a^{[q]_{d-2}}c_{d-1}F^{d-2} + a^{[q]_{d-3}}c_{d-1}c_{d-2}F^{d-3} + ... + a^{[q]_0}c_{d-1}c_{d-2}...c_1.$$

This is verified by checking that the coefficients of $F^n$ of each side of the identity agree for all $n$. This ugly calculation is omitted. (In fact, the identity was discovered by extensive computer algebra calculations with *Mathematica* of small degree cases.) We apply both sides of the identity to $x$. This yields $a^{[q]_d}v^{q^d} - c_0c_1c_2...c_{d-1}v = 0$. Hence $v^{q^d} = ((c_0c_1c_2...c_{d-1})/a^{[q]_d})v$, and we are done, if we can show that $L' = K'(v)$ (note that this will also show that $v \neq 0$). This follows from the following lemma.

**Lemma** The (right) greatest common divisor of $\phi(aF+b)$ and $\psi(aF+b)$ is 1, i.e. they are (right) relatively prime.

*Proof* As described in 1.10 of [G], the greatest common divisor is calculated as follows. Let $W_\phi$ and $W_\psi$ denote the set of zeros in $\overline{K}$ of $(\phi(aF+b))x = 0$ and $(\psi(aF+b))x = 0$ respectively. If $W = W_\phi \cap W_\psi$, then the greatest common divisor is the additive polynomial $\prod_{\alpha \in W}(x-\alpha)$. We therefore need to show that $W = \{0\}$. This is accomplished by using the easily verified identity

$$\phi(aF+b) = -\zeta\psi(aF+b) + \psi(aF+b)(aF+b).$$

Suppose that $u \in W, u \neq 0$. By the last identity, $(\psi(aF+b))(aF+b)u = 0$. Since the coefficients of $\phi$ are in $\mathbf{F}_q$, $(\phi(aF+b))(aF+b)u = 0$, so $aF+b$ is an endomorphism of $W$, i.e. $W$ is an $\mathbf{F}_q[aF+b]$-submodule of $W_\phi$. Since $W_\phi$ is 1-dimensional over $\mathbf{F}_q[aF+b] = \mathbf{F}_{q^d}$, $W = W_\phi$, which contradicts the fact that $\#W_\psi < \#W_\phi$.

This incidentally shows that the identity in the first lemma above in fact gives the least common multiple of $\phi(aF+b)$ and $\psi(aF+b)$ - see section 1.10 of [G], where consequences of the existence of a right division algorithm in Ore rings are discussed.

By the lemma, we can find polynomials $j(F), k(F) \in K'\{F\}$ such that $j(F)\psi(aF+b) + k(F)\phi(aF+b) = 1$. Applying this to $x$ gives $j(F)v = x$, so $x \in K'(v)$ and since $v = (\psi(aF+b))x, v \in K'(x)$.

## 3. The cases $d = 1$ and $d = 2$

The above lemma allows us to show that every representation $G_K \to GL_1(\mathbf{F}_{q^d})$ is Drinfeld if $d = 1$ or 2, except for one special case for $d = 2$, namely when $K = \mathbf{F}_q$ and the image of the representation is in $GL_1(\mathbf{F}_q)$. The idea is to let $L$ be the fixed field of its kernel and to show that $L = L_{a,b,\phi}$ for some $a, b \in K$ and irreducible $\phi \in \mathbf{F}_q[T]$ of degree $d$. Note that this is enough to show that the associated representation is Drinfeld since the property of being Drinfeld depends only on the field cut out - the representation can be changed by picking a different basis for the corresponding $V$.

(I) $d = 1$.

Given representation $G_K \to GL_1(\mathbf{F}_q)$, we let $L$ be the fixed field of its kernel. Then $L/K$ is a Kummer extension and so is of the form $L = K(v)$, where $v^{q-1} = c \in K$.

Taking $a = 1, b = -c$, and $\phi(T) = T$ (so that $\zeta = 0$), we get by the Drinfeld construction a representation that, by the last lemma, yields $L_{a,b,\phi} = L$ (since $(\zeta - b)/a = c$).

(II) $d = 2$.

There are now two cases, namely according as $\zeta \in K$ or $\zeta \notin K$.

Case (i): $\zeta \in K$.

Then $\mathbf{F}_q(\zeta) = \mathbf{F}_{q^d} \leq K$ and so $L/K$ is a Kummer extension, say $L = K(v)$ with $v^{q^2-1} = c \in K$. We wish to find $a, b \in K$ such that $(\zeta - b)(\zeta - b^q)/a^{q+1} = c$. Note that $(\zeta - b)(\zeta - b^q)/a^{q+1} = \frac{\zeta - b}{\zeta^q - b}(\frac{\zeta^q - b}{a})^{q+1}$, so if we set $b = (c\zeta^q - \zeta)/(c - 1)$ and $a = \zeta^q - b$, then this all simplifies to $c$. We just have to make sure that $c \neq 1$, but $c$ is only defined up to a $(q^2 - 1)$th power, so we have the necessary flexibility, unless $K = \mathbf{F}_{q^2} = L$. In that case, we need to pick $b \in K$ such that $(\zeta - b)(\zeta - b^q)$ is a $(q+1)$th power in $K^*$, i.e. is a nonzero element of $\mathbf{F}_q$. This is accomplished in exactly the same way as described in the first paragraph of I below.

Case (ii): $\zeta \notin K$.

The idea is to show that the process, considered in Lemma 1, of obtaining $L'$ as the compositum of $K' = K(\zeta)$ and $L$ can be suitably reversed.

We begin with the cases when $L$ and $K(\zeta)$ are disjoint, i.e. $\zeta \notin L$. The extension $L'/K$ is Galois with Galois group $<\sigma> \times <\tau>$ where $\sigma$ has order 2 and $\tau$ has order $m$ dividing $q^2 - 1$. The fixed fields of $\sigma$ and $\tau$ are $L$ and $K' = K(\zeta)$ respectively.

The extension $L'/K'$ is Kummer and so $L' = K'(v)$ for some $v$ such that $v^{q^2-1} = c \in K'$. We claim that there exist $a, b \in K$ such that $((\zeta - b)(\zeta - b^q))/a^{q+1} = c$. The argument goes as follows.

Let $w = \sigma(v)$. Then $w^{q^2-1} = \sigma(c)$. Suppose, without loss of generality, that $\tau(v) = \eta v$, where $\eta$ is an $m$th root of unity in $K'$. The fact that $\sigma$ and $\tau$ commute, implies that $\tau(w) = \eta^q w$. Let $y = wv^{-q}$. We check that $\tau(y) = y$ and so $y \in K'$. We calculate that $\sigma(y)y^q c = 1$.

At this point, we have a division into two cases depending on whether $y \in K$ or not.

I. Say $y \in K$. Then $\sigma(y) = y$. Hence, $c = (1/y)^{q+1}$ is the $(q+1)$th power of an element of $K$ and so, to write $c$ in the form $(\zeta - b)(\zeta - b^q)/a^{q+1}$ (up to $(q^2 - 1)$th powers of elements of $K'$), we must equivalently be able to pick $b \in K$ such that $(\zeta - b)(\zeta - b^q)$ is a $(q+1)$th power in $K$ times a $(q^2 - 1)$th power in $K'$. This can be done so long as $K \neq \mathbf{F}_q$. For instance, in the case of odd characteristic, suppose $\phi = T^2 - \lambda$. Pick any $u \in K - \mathbf{F}_q$. Set $b = (u^{q+1} - \lambda)/(u^q - u)$ Then $(\zeta - b)(\zeta - b^q)/a^{q+1} = ((\zeta - u)(\zeta - u^q)/((u^q - u)a))^{q+1} = ((u - b)/a)^{q+1}(\zeta(u + \zeta))^{q^2 - 1}$, which is of the desired form.

If $K = \mathbf{F}_q$, then since $c$ is a $(q+1)$th power of an element $1/y$ of $K$, we can pick $v$ so that $v^{q-1} = 1/y \in K$. Then $L = K(v)$ has degree dividing $q - 1$ over $K$. Suppose now $(\zeta - b)(\zeta - b^q) = k^{q+1}r^{q^2-1}$ for some $b, k \in K, r \in K'$. Since $K' = \mathbf{F}_{q^2}$, $r^{q^2-1} = 1$. Moreover, $k^{q+1} = k^2$ and $b^q = b$ since they are in $K$. The equation reduces to $(\zeta - b)^2 = k^2$, so $\zeta - b = \pm k$, impossible since $\zeta \notin K$.

II. Say $y \notin K$. Since $1/\sigma(y) \in K' - K$ and $K' = K(\zeta)$ has degree 2 over $K$, we can

write $1/\sigma(y) = s\zeta - r$ with $r, s \in K, s \neq 0$. Then $(s\zeta - r)(s^q\zeta - r^q) = 1/(\sigma(y)y^q) = c$. Let $b = r/s$ and $a = 1/s$. We have shown that $((\zeta - b)(\zeta - b^q))/a^{q+1} = c$.

It follows that $L'$ is the compositum of $L_{a,b,\phi}$ and $K'$. The fixed field of $\sigma$ equals $L$ and $L_{a,b,\phi}$ and so the two fields must coincide.

It remains to deal with the case when $L$ and $K(\zeta)$ are not disjoint, i.e. $\zeta \in L - K$. In this case we have a tower of fields with, say, $\text{Gal}(L/K) = <\sigma>$ so that $\text{Gal}(L/K(\zeta)) = <\sigma^2>$. Since $L/K(\zeta)$ is Kummer, there is $v$ such that $\sigma^2(v) = \eta v$ with $\eta$ an $m$th root of unity, where $m = [L : K(\zeta)]$. Note that since $[L : K] = 2m$ divides $q^2 - 1$, $\eta$ is a square in $\mathbf{F}^*_{q^2}$. We can write $\eta = \mu^{2q}$ then with $\mu \in \mathbf{F}^*_{q^2}$.

Setting $y = v^q \sigma(v)^{-1} \mu$, we check that $\sigma(y)y^q = v^{q^2-1} = c$, say. So long as $y \notin K$, we can pick $a, b \in K$ such that $(\zeta - b)/a = \sigma(y)$ and we are done. The case of $y \in K$ is handled exactly as in I above.

**Corollary** Every cyclic extension of $K \neq \mathbf{F}_q$ of degree dividing $q^2 - 1$ is the splitting field of an equation of the form

$$a^{q+1}x^{q^2-1} + a(b^q + b + 1)x^{q-1} + (b^2 + b + \lambda) \quad (\text{char}(K) = 2)$$

$$a^{q+1}x^{q^2-1} + a(b^q + b)x^{q-1} + (b^2 - \lambda) \quad (\text{char}(K) \neq 2),$$

where $\lambda \in \mathbf{F}_q$ is chosen so that $T^2 + T + \lambda$, respectively $T^2 - \lambda$, is irreducible over $\mathbf{F}_q$.

*Proof* In the case of characteristic two, pick $\lambda \in \mathbf{F}_q$ such that $\phi(T) = T^2 + T + \lambda$ is irreducible over $\mathbf{F}_q$. Then $\phi(aF+b) = a^{q+1}F^2 + a(b^q+b+1)F + (b^2+b+\lambda)$. Applying this to $x$ and dividing by $x$ yields the desired equation. The odd characteristic case proceeds similarly with $\phi(T) = T^2 - \lambda$ ($\lambda$ chosen to make $\phi$ irreducible over $\mathbf{F}_q$).

Note that the corollary still holds if $K = \mathbf{F}_q$ and the degree does not divide $q-1$. If the degree does divide $q-1$, then the extension is Kummer and so a splitting field for e.g. $ax^{q-1} + b$.

Partway through the main result of this section, the choice $b = (u^{q+1} - \lambda)/(u^q - u)$ was made and it is probably worthwhile explaining where this choice comes from. We explain this in the case of odd characteristic. A similar approach works in even characteristic.

Setting $y = ax^{q-1}$ in the second equation of the corollary leads to $y^{q+1} + (b^q + b)y + (b^2 - \lambda) = 0$. The field $K(y)$ is an important intermediate field between $L_{a,b,\phi}$ and $K$, as evidenced in the next section. The choice of $b$ that we are discussing, is one that will ensure that the equation for $y$ splits completely. The idea is to set $y = u - b$, which yields $(u^{q+1} - \lambda) - b(u^q - u) = 0$ (†). Hence the choice of $b$. The fact that the equation splits completely can be forcefully seen by the next lemma.

**Lemma** Let $\mu = 1/\lambda$ and $H_\mu$ be the image in $PGL_2(\mathbf{F}_q)$ of the nonsplit Cartan subgroup

$$\left\{ \begin{pmatrix} \alpha & \beta \\ \mu\beta & \alpha \end{pmatrix} : (\alpha, \beta) \neq (0,0) \right\},$$

a cyclic group of order $q+1$. Then

$$(u^{q+1} - \lambda) - b(u^q - u) = \prod_{\sigma \in H_\mu} (u - \sigma(U)),$$

where $U$ is one root of (†) and $\sigma$ acts by fractional linear transformation.

*Proof* The proof follows automatically by checking that $\sigma(U)$ satisfies (†).

## 4. Genus Constraints

In this section, we consider properties of $L_{a,b,\phi}/K$. Without loss of generality, we can replace $K$ by its subfield $\mathbf{F}_q(a,b)$. We then have $L = K(x)$, where $x$ satisfies $(\phi(aF + b))x = 0$. Setting $y = ax^{q-1}$, we get an intermediate field $M = K(y)$. The important facts are that the extension $L/M$ is Kummer and that the equation satisfied by $y$ has coefficients involving $b$ but not $a$ (see e.g. the equations in the above corollary). It is therefore sufficient for our purposes to study the case of $K = \mathbf{F}_q(b)$ and $\phi = T + b$. (This case of the Drinfeld construction was first considered by Carlitz and, in greater detail, by Hayes [H].)

The idea is to calculate the genus of a certain field $N$ of degree $(q^d-1)/(k(q-1))$ over $\mathbf{F}_q(b)$.

**Lemma** The genus of $N$ is

$$g_N = \frac{1}{2}(d-2)((q^d-1)/(k(q-1)) - 1).$$

*Proof* By Riemann-Hurwitz, $2g_N - 2 = -2(\frac{q^d-1}{k(q-1)}) + \deg(\mathcal{D})$, where $\mathcal{D} = \mathcal{B}^s$, where $\mathcal{B}$ is the totally ramified prime of $N$ over $(\phi)$, and $s = e - 1 = \frac{q^d-1}{k(q-1)} - 1$. Then $\deg(\mathcal{B}) = d$ implies that $2g_N - 2 = -2(\frac{q^d-1}{k(q-1)}) + d(\frac{q^d-1}{k(q-1)} - 1)$, whence the result.

**Corollary** $g_N = 0$ if and only if $d = 1$ or $d = 2$ or $k = \frac{q^d-1}{q-1}$.

**Theorem** If $\rho: G_K \to GL_1(\mathbf{F}_{q^d})$ is Drinfeld, then (i) $d = 1$ or (ii) $d = 2$ or (iii) $\pi \circ \rho$ surjects onto $GL_1(\mathbf{F}_{q^d})/GL_1(\mathbf{F}_q)$, where $\pi$ is the quotient map $GL_1(\mathbf{F}_{q^d}) \to GL_1(\mathbf{F}_{q^d})/GL_1(\mathbf{F}_q)$.

Take $k$ to be $\#(\pi \circ \rho(G_K))$. Then $b$ is such that $N$ specializes to $K$ and so by Lüroth, $g_N = 0$, leading to the desired result.

There are two important consequences to this, first that Drinfeld representations tend to have large images (results like this were already established by Goss [G], section 7.7) and second that representations that are not Drinfeld certainly exist (by picking $d > 2$ and taking a representation which does not surject onto $GL_1(\mathbf{F}_{q^d})/GL_1(\mathbf{F}_q)$). In the next section, we show that there are many representations that are not Drinfeld but that are surjective.

## 5. Surjective Representations That Are Not Drinfeld

Take $q = 2, d = 3$, and $\phi = T^3 + T + 1$. We assume that $\phi$ is irreducible over $K$ - in other words any root $\zeta$ of $\phi$ satisfies $\zeta \notin K$. Then $\text{Gal}(K(\zeta)/K) = <\sigma>$ has order 3. We will always make the choice of $\sigma$ such that $\sigma(\zeta) = \zeta^2$ (be warned that for fields $\mathbf{F}_{2^m}$ with $m \equiv 2 \pmod{3}$ this makes $\sigma$ the square of Frobenius). We provide a method (that in fact generalizes to any $d > 2$ and to other $q$) of obtaining numerous representations that are not Drinfeld, so long as $K$ does not satisfy a certain hypothesis (A) below.

**Hypothesis (A)** Let $S = \{\sigma(x)x^{-2} : x \in K(\zeta)^*\}$, a subgroup of the multiplicative group of $K(\zeta)$. Say that $K$ satisfies hypothesis (A) if every coset of $S$ in $K(\zeta)^*$ contains an element of the form $r + s\zeta (r, s \in K)$.

**Theorem** Suppose that $K$ does not satisfy hypothesis (A). Then there exists a **surjective** representation (in fact many such) $G_K \to GL_1(\mathbf{F}_{q^d})$, that is not Drinfeld.

*Proof* Let $f(x) = x\sigma^{-1}(x^2)\sigma^{-2}(x^4)$, a homomorphism of the multiplicative group of $K' = K(\zeta)$ to itself. Note that $f$ satisfies two useful identities, (i) $\sigma(f(x)) = f(x)^2 \sigma(x)^{-7}$ and (ii) $x^7 = f(x)^{-1}\sigma^{-1}(f(x))^2$.

Pick $y \in K'$ such that the coset of $y$ contains no element of the form $r+s\zeta (r, s \in K)$. Let $c = f(y)$ and $L' = K(v)$ with $v^7 = c$. Then $L'/K$ is Galois with Galois group $C_3 \times C_7$. (Note that $v \notin K$, since otherwise $f(v) = v^7 = f(y)$ and, by the injectivity of $f$ proven below, $y = v \in K$, a contradiction.)

We claim that $c$ is not of the form $((\zeta - b)(\zeta - b^2)(\zeta - b^4))/a^7$ times a 7th power of an element of $K'$ for any $a, b \in K$, and so the subfield $L$ of degree 7 over $K$ is not obtained by the Drinfeld construction and we are done.

We first show that $f$ is injective. Suppose that $x \in K'$ satisfies $f(x) = 1$. By identity (ii), we get that $x^7 = 1$ and so $x = \zeta^i$ for some $i$. Since $f(\zeta^i) = \zeta^{3i}$, it follows that $x = \zeta^i = 1$.

If the subfield $L$ of degree 7 over $K$ is obtained by the Drinfeld construction, then $c$ is of the form $k^7((\zeta - b)(\zeta - b^2)(\zeta - b^4))/a^7$ for some $a, b \in K, k \in K'$. We check that $x = k^{-1}\sigma^{-1}(k^2)$ is a solution of $f(x) = k^7$ and so, by the injectivity of $f$, is the unique such solution. Then, $f(k^{-1}\sigma^{-1}(k^2)(\zeta - b)/a) = c = f(y)$, and so by the injectivity of $f$, $(\zeta - b)/a = yk\sigma^{-1}(k^{-2})$, which contradicts our choice of $y$.

It remains to make some comments on what fields $K$ satisfy hypothesis (A) and what fields do not.

**Lemma** If the equation

$$(Q) \quad u^4 + u^3v + u^2v^2 + uv^3 + v^4 + u^3w + u^2vw + v^3w + v^2w^2 + vw^3 + w^4 = 0$$

has no solutions in $K$ other than $(0, 0, 0)$, then $K$ fails to satisfy hypothesis (A).

*Proof* Suppose that $K$ satisfies (A). Then $\zeta^2$ is in the same coset of $S$ as some $r + s\zeta (r, s \in K)$. So there is some $x = u + v\zeta + w\zeta^2 (u, v, w \in K$ not all 0) such that $\zeta^2 = \sigma(x)x^{-2}(r + s\zeta)$. Writing this in terms of $u, v, w$ and clearing denominators, we get, by comparing coefficients of $1, \zeta, \zeta^2$, three linear equations in $r, s$. We use two of these to solve for $r, s$ and plug in the third to get that some expression in $u, v, w$ is 0. The numerator of that expression is (Q). This provides the desired contradiction.

This lemma is very useful in establishing that certain fields fail to satisfy (A). With a little more work, we can establish a converse. As in the above proof, we might ask whether $a + b\zeta + c\zeta^2$ is in the same coset as some $r + s\zeta (a, b, c, r, s \in K)$. Proceeding as above yields a quartic $Q(a, b, c) : cu^4 + au^3v + cu^3v + bu^2v^2 + cu^2v^2 + auv^3 + buv^3 + cuv^3 + cv^4 + au^3w + bu^3w + cu^3w + cu^2vw + auv^2w + av^3w + cv^3w + au^2w^2 + auvw^2 + bv^2w^2 + cv^2w^2 + buw^3 + avw^3 + bvw^3 + cvw^3 + cw^4 = 0$, with $Q(0, 0, 1)$ being $Q$ above.

**Lemma** If $Q(a, b, c)$ has no solutions in $K$ other than $(0, 0, 0)$ for some choice of $a, b, c \in K$ with $c \neq 0$, then $K$ fails to satisfy hypothesis (A). If $K$ fails to satisfy hypothesis (A), then there is some choice of $a, b, c \in K$ for which $Q(a, b, c)$ has no nontrivial solutions.

*Proof* Exactly as for the previous lemma.

**Lemma** If $c \neq 0$, then $Q(a, b, c)$ defines a nonsingular curve over $K$.

*Proof* We first calculate the partial derivatives of $Q(a, b, c)$ with respect to $u, v, w$ and set these equal to 0. These equations are linearly dependent and we solve them to obtain $a, b$ in terms of $c$. Without loss of generality we may assume that $c = 1$. Plugging the expressions for $a, b$ into $Q(a, b, c)$ yields a homogeneous nonic, namely $u^9 + u^7v^2 + u^6v^3 + u^5v^4 + u^2v^7 + uv^8 + v^9 + u^7vw + u^3v^5w + uv^7w + u^7w^2 + u^6vw^2 + u^5v^2w^2 + uv^6w^2 + v^7w^2 + u^6w^3 + u^3v^3w^3 + u^2v^4w^3 + v^6w^3 + u^5w^4 + u^3v^2w^4 + uv^4w^4 +$

$v^5w^4+u^3vw^5+u^2v^2w^5+uv^3w^5+u^2vw^6+u^2w^7+uvw^7+v^2w^7+uw^8+vw^8+w^9 = 0$.
This nonic actually equals $(u^3 + uv^2 + v^3 + uvw + uw^2 + vw^2 + w^3)^3 = 0$, which factors as $((u + \zeta v + \zeta^2 w)(u + \zeta^2 v + \zeta^4 w)(u + \zeta^4 v + \zeta w))^3 = 0$. These are lines defined only over fields containing $\zeta$, so not $K$.

**Theorem** The only finite fields of characteristic 2 that fail to satisfy (A) are $\mathbf{F}_2$ and $\mathbf{F}_4$.

*Proof* The homogeneous quartics $Q(a, b, c)$ have genus $g \leq 3$ (in fact exactly 3) as smooth projective curves over $\mathbf{F}_2$. Letting $N(q)$ denote the number of points on the curve over $\mathbf{F}_q$, by the Hasse-Weil bound $N(q) \geq q + 1 - 2g\sqrt{q}$, which is $> 0$ if $q \geq 36$. We check $q = 2, 4, 16, 32$ separately (note that $q \neq 8$ by the assumption of irreducibility of $\phi$). On the one hand, for quartic $Q$, $N(2) = 0$ and $N(4) = 0$, and by an above lemma $\mathbf{F}_2$ and $\mathbf{F}_4$ do not satisfy (A). Using *Magma,* it is easy to check that $\mathbf{F}_{16}$ and $\mathbf{F}_{32}$ satisfy (A).

**Theorem** Let $K$ be a field of characteristic 2 with a discrete valuation such that its residue field $k$ fails to satisfy (A). Then $K$ fails to satisfy (A).

*Proof* Suppose that $K$ satisfies (A). Then $Q(a, b, c)$ has a solution $u, v, w \in K$ (whatever $a, b, c \in K$ are, $c \neq 0$). Then their leading terms (in $k$) also yield a nonzero solution of $Q(a, b, c)$. A contradiction.

It follows then that if $K$ fails (A), then so do $K(T), K(T_1, ..., T_n), K((T)), ...$, as well as function fields with constant field $\mathbf{F}_2$ with at least one place of degree 1 or 2, and so on.

## 6. Higher degree representations

Cases where $r > 1$ are poorly understood, except in one instance, namely when the given representation is into $GL_r(\mathbf{F}_q)$. In that case, we can say the following.

**Theorem** Let $K$ be infinite and $\rho : G_K \to GL_r(\mathbf{F}_q)$ be a representation. Then $\overline{\rho}$ is Drinfeld. This is not necessarily true if $K$ is finite.

*Proof* Suppose that $K$ is infinite. Let $L$ be the fixed field of the kernel of $\rho$. Let $H$ denote $\text{Gal}(L/K)$, which is isomorphic to the image of $\rho$. Let $V$ be the $\mathbf{F}_q[H]$-module corresponding to the embedding of $H$ in $GL_r(\mathbf{F}_q)$. By the normal basis theorem, $V$ embeds $\mathbf{F}_q[H]$-linearly in the additive group $L^+$ of $L$ (since it contains free $\mathbf{F}_q[H]$-modules of arbitrarily high finite rank and by duality for group rings these are also cofree of arbitrary finite rank). Let $g(x) = \prod_{\alpha \in V}(x - \alpha)$. Since $V$ is an $\mathbf{F}_q$-vector space, the polynomial $g$ is indeed additive and so lies in $K\{F\}$.

Setting $\phi(T) = T$, the $T$-division points are the roots of $g$, i.e. $V$, with the given action. Finally, the extension of $K$ generated by the elements of $V$, $K(V)$, is indeed $L$, since $\rho$ factors through $\text{Gal}(K(V)/K)$.

Suppose that $K$ is finite. If $\rho$ is Drinfeld, then $\phi$ has degree $d = 1$, say $\phi = aT+b$. Then $\phi(g(F)) = ag(F) + b = h(F)$, say, so $V$ is the set of zeros in $\overline{K}$ of $h(F)x = 0$ and is an $\mathbf{F}_q$-subspace of $L^+$, where $L$ is the fixed field of the kernel of $\rho$. The action of $H = \text{Gal}(L/K)$ on $L^+$ restricts to $V$ to produce $\rho$, but for large $r$, $V$ will not embed in $L^+$, which is a free $\mathbf{F}_q[H]$-module of rank $[K : \mathbf{F}_q]$.